\documentclass{amsart}
\usepackage{graphicx}
\usepackage{latexsym}
\usepackage{hyperref}
\usepackage{amsfonts}
\usepackage[all]{xy}
\usepackage{amssymb, mathrsfs, amsfonts, amsmath}
\usepackage{amsbsy}
\usepackage{amsfonts}
\newtheorem*{acknowledgement}{Acknowledgement}
\newtheorem{corollary}{Corollary}

\newtheorem{theorem}{Theorem}

\numberwithin{equation}{section}

\begin{document}
\title[Einstein warped products]{Bounds on volume growth of geodesic balls\\ for Einstein warped products}
\author{A. Barros$^{1}$, R. Batista$^{2}$ \& E. Ribeiro Jr.$^{3}$}
\address{$^{1}$ Universidade Federal do Cear\'a - UFC, Departamento  de Matem\'atica, Campus do Pici, Av. Humberto Monte, Bloco 914, 60455-760-Fortaleza /CE , Brazil} \email{abbarros@mat.ufc.br} 

\address{$^{2}$ Universidade Federal do Piau\'i - UFPI, Departamento de Matem\'{a}tica, Campus Petr\^onio Portella, 64049-550-Teresina /PI, Brazil} \email{rmarcolino@ufpi.edu.br}
\thanks{$^{1,2}$ Partially supported by grant from CNPq/Brazil}

\address{$^{3}$ Current: Department of Mathematics, Lehigh University, Bethlehem, PA, 18015, United States\\ Permanent: Universidade Federal do Cear\'a - UFC, Departamento  de Matem\'atica, Campus do Pici, Av. Humberto Monte, Bloco 914, 60455-760-Fortaleza /CE, Brazil}\email{ernani@mat.ufc.br}
\thanks{$^{3}$ Partially supported by grants from  PJP-FUNCAP/Brazil and CNPq/Brazil}

\date{May 2, 2014}
\begin{abstract}
The purpose of this note is to provide some volume estimates for Einstein warped products similar to a classical result due to Calabi and Yau for complete Riemannian manifolds with nonnegative Ricci curvature. To do so, we make use of the a\-pproach of quasi-Einstein manifolds which is directly related to Einstein warped pro\-duct. In particular,  we present an obstruction for the existence of such a class of manifolds.
\end{abstract}

\maketitle

\section{Introduction}
It has long been a goal of mathematicians to understand the geometry of Einstein ma\-ni\-folds as well as Einstein-type manifolds, for instance Ricci solitons and
$m$-quasi-Einstein manifolds. Surely, this is a fruitful problem in Riemannian Geometry. Ricci solitons model the formation of
singularities in the Ricci flow and  they correspond to self-similar solutions for this flow; for more details in
this subject we recommend  the survey by Cao \cite{Cao}. On the other hand, one of the motivation to
study $m$-quasi-Einstein metrics on a Riemannian manifold is its direct relation with Einstein warped product. For comprehensive references on such a theory, see \cite{brG}, \cite{besse}, \cite{Case2}, \cite{Corvino}, \cite{petersen}, \cite{Kim} and \cite{WW}.

One fundamental ingredient to understand the behavior of Einstein warped product is the $m$-Bakry-Emery Ricci tensor which appeared
previously in \cite{bakry} and \cite{qian} as a modification of the classical Bakry-Emery tensor $Ric_{f}=Ric+\nabla^2 f.$ More exactly, the $m$-Bakry-Emery Ricci tensor is given by
\begin{equation}
\label{bertens}
Ric_{f}^{m}=Ric+\nabla ^2f-\frac{1}{m}df\otimes df,
\end{equation}where $f$ is a smooth function on $M^n$ and $\nabla ^2f$ stands for the Hessian form.

A Riemannian manifold $(M^n,\,g),$ $n\geq 2,$ will be called $m$-{\it quasi-Einstein manifold}, or simply {\it quasi-Einstein manifold}, if there exist a smooth potential function $f$
on $M^n$ and a constant $\lambda$ satisfying the following fundamental equation
\begin{equation}
\label{eqqem}
Ric_{f}^{m}=Ric+\nabla ^2f-\frac{1}{m}df\otimes
df=\lambda g.
\end{equation}
When $m$ goes to infinity, equation
(\ref{eqqem}) reduces to the one associated to a Ricci soliton.
Furthermore, when $m$ is a positive integer it
corresponds to a warped product Einstein metric, for more details see,
for instance, \cite{Case2}.  Following the terminology of Ricci
solitons, a  $m$-quasi-Einstein metric $g$ on a manifold $M^n$ will be
called \emph{expanding}, \emph{steady} or \emph{shrinking},
respectively, if  $\lambda<0,\,\lambda=0$ or $\lambda>0$. Moreover,
a $m$-quasi-Einstein manifold will be called \emph{trivial} if $f$ is constant, otherwise it will be \emph{nontrivial}. We notice that
the triviality implies that $M^n$ is an Einstein manifold.

In order to proceed we recall a classical result which gives a characterization for Einstein warped products.

\begin{theorem}[\cite{besse}, \cite{Kim}]
\label{1teo}
Let $M^{n}\times_{u}F^{m}$ be an Einstein warped product with Einstein constant $\lambda$, warping function $u=e^{-\frac{f}{m}}$ and Einstein fiber $F^{m}.$ Then the weighted manifold $(M^{n}, g_{M}, e^{-f} dM)$ satisfies the $m$-quasi-Einstein equation (\ref{eqqem}). Furthermore the Einstein constant $\mu$ of the fiber satisfies
\begin{equation}\label{2eq}
\Delta f-|\nabla f|^{2} = m\lambda-m\mu e^{\frac{2f}{m}}.
\end{equation}
Conversely, if the weighted manifold $(M^{n}, g_{M}, e^{-f}dM)$ satisfies (\ref{eqqem}), then $f$ satisfies (\ref{2eq}) for some constant $\mu\in\mathbb{R}$. Considering the warped product $N^{n+m}=M^{n}\times_{u}F^{m}$, with $u=e^{-\frac{f}{m}}$ and Einstein fiber $F$ with $Ric_{F}=\mu g_{F}$, then $N$ is also Einstein with $Ric_{N}=\lambda g_{N}$, where $g_{N}=g_{M}+u^2g_{F}$.
\end{theorem}

Clearly, Theorem \ref{1teo} shows that a $m$-quasi-Einstein structure provides a structure of Einstein warped product; for more details we recommend \cite{Kim}. Therefore, classifying $m$-quasi-Einstein manifolds or understanding their geometry is definitely an important issue.

One should point out that some examples of expanding $m$-quasi-Einstein manifolds with arbitrary $\mu$ as well as steady quasi-Einstein manifolds with $\mu>0$ were constructed in \cite{besse}. While Case  \cite{Case} showed that steady $m$-quasi-Einstein manifolds with $\mu\leq0$ are trivial. In \cite{qian},  Qian proved that shrinking $m$-quasi-Einstein manifolds must be compact. Moreover, as already noticed in \cite{Kim} the converse result is also true. From what it follows that a $m$-quasi-Einstein manifold is compact if and only if $\lambda>0.$ See also \cite{WW} for further discussion.

Based on the above results and inspired by works due to Calabi and Yau, we shall investigate bounds on volume growth of geodesic balls for non compact $m$-quasi-Einstein manifolds, in particular $\lambda$ must be non positive. For our purposes we recall that Calabi \cite{calabi} and Yau  \cite{yau} proved that every metric with non negative Ricci tensor on a non compact smooth manifold satisfies
\begin{equation}
\label{aim}
V ol(B_{p}(r))\geq c r
\end{equation} for any $r>r_{0}$ where $r_{0}$ is a positive constant and $B_{p}(r)$ is the geodesic ball of radius $r$ centered at $p\in M^n$ and, $c$ is a constant does not depend on $r.$ In a similar way Munteanu and Sesum \cite{Mun} obtained the same type of growth for steady gradient Ricci soliton.

Now we state our first result concerning the growth of volume of geodesic balls  for non compact $m$-quasi-Einstein manifolds which is similar to Calabi-Yau estimate.

\begin{theorem}\label{thmS}
Let $\big(M^{n},\,g,\,f\big)$ be a non compact steady $m$-quasi-Einstein manifold with $m \in(1,\infty].$ Then there exist constants $c$ and $r_{0}>0$ such that for any $r>r_{0}$
\begin{equation}\label{eq6}
 Vol(B_{p}(r))\geq cr.
\end{equation}
\end{theorem}

This immediately yields the following corollary.

\begin{corollary}
\label{cor1thmS}
Let $N=M^{n}\times_{u} F^{m}$ be a Ricci flat warped product. Then there exist constants $c$ and $r_{0}>0$ such that for any $r>r_{0}$
\begin{equation}\label{eq61}
Vol(B_{p}(r))\geq cr,
\end{equation}for geodesic balls of the base.
\end{corollary}

One question that naturally arises from the above results is to know what occurs on expanding $m$-quasi-Einstein manifold. In this case, we obtain the following volume growth of geodesic balls.

\begin{theorem}\label{thmE}
Let $\big(M^{n},\,g,\,f,\lambda\big)$ be a non compact expanding $m$-quasi-Einstein manifold with $m \in(1,\infty)$ and  $\mu\leq 0.$ Supposing that $f \geq -k$ for some positive constant $k,$ then there exist constants $c$ and $r_{0}>0$ such that for any $r>r_{0}$
\begin{equation}\label{eq7}
 Vol(B_{p}(r))\geq cr.
\end{equation}
\end{theorem}

In the sequel, inspired by ideas developed in \cite{detang},  \cite{Li} and \cite{morgan} we shall prove an $f$-volume estimate of geodesic balls on expanding $m$-quasi-Einstein manifolds. More precisely, we have the following result.

\begin{theorem}\label{thmE-C}
Let $\big(M^{n},\,g,\,f,\lambda\big)$ be a non compact expanding $m$-quasi-Einstein manifold with $m \in[1,\infty)$ and  $\mu=0.$ Then there exists a constant $c$ such that for any $r>1$
\begin{equation}\label{eq8}
  Vol_{f}(B_{p}(r))\geq c e^{\sqrt{-\lambda m}r}.
\end{equation}
\end{theorem}

It should be emphasized that there are  several further interesting obstructions to the existence of Einstein metrics. Based in this, in \cite{besse} (cf. page 265) it was posed the following question:

\begin{flushright}
\begin{minipage}[t]{4.37in}
 \emph{``Does there exist a compact Einstein warped product with non constant warping
function?"}
 \end{minipage}
\end{flushright}

In 2003 Kim and Kim \cite{Kim} gave a partial answer for this question by means of the quasi-Einstein approach. While  Case \cite{Case} studied this problem without compactness assumption for $m$-quasi-Einstein manifolds with $\lambda=0$ and $\mu\leq0.$ Here, we use the weak Maximum Principle at infinity for the $f$-Laplacian to obtain a triviality result for Einstein warped products with negative scalar curvature. More precisely, we have the following result.

\begin{theorem}\label{thmPW}
Let $N = M^{n}\times_{u} F^{m}$ be a complete Einstein warped product with Einstein constant $\lambda<0,$ warping function $u$ and Einstein constant of the fiber $F$ satisfying $\mu<0.$ If the warping function satisfies $$u\leq\sqrt{\frac{2\mu}{\lambda}},$$
then $u$ is a constant function and $N$ is a Riemannian product.
\end{theorem}

\section{Proof of the results}

In order to set the stage for the proofs to follow let us recall some classical equations. First, considering the function $u=e^{-\frac{f}m}$ on $M^n,$ we immediately have $\nabla u=-\frac{u}{m}\nabla f$ as well as
\begin{equation}\label{eq5}
\nabla^{2}f-\frac{1}{m}df\otimes df=-\frac{m}{u}\nabla^{2}u.
\end{equation}Next, a straightforward computation involving (\ref{eqqem}) and (\ref{2eq}) gives
\begin{equation}
\frac{u^{2}}{m}(R-\lambda n)+(m-1)|\nabla u|^{2}=-\lambda u^{2}+\mu.
\end{equation}
We also recall that Wang \cite{wang} proved that if $\lambda\leq 0,$ then $R\geq \lambda n.$ From what it follows that
\begin{equation}\label{eq8}
(m-1)|\nabla u|^{2}\leq-\lambda u^{2}+\mu.
\end{equation}

Now we are ready to prove the results.

\subsection{Proof of Theorem \ref{thmS}}

\begin{proof}
To begin with, we notice that when $m=\infty$ we have a gradient Ricci soliton and in this case the result follows from \cite{Mun}. From now on we can assume that $m\in(1,\infty).$ From this, since $\lambda=0$ we get
\begin{equation}\label{eqvol4}
|\nabla u|^{2}\leq\frac{\mu}{m-1}.
\end{equation} Taking into account that $R\geq 0$ and $u>0$ we deduce $$\int_{B_{p}(r)} uRd\sigma \geq 0$$ for each $r>0,$  where $d\sigma$ denotes the Riemannian volume form. Consequently, if for all $r>0$ we have $$\int_{B_{p}(r)} uRd\sigma=0,$$ then $R=0$ on $M^n.$

On the other hand, Wang  \cite{wang} (see also \cite{Case2}) proved that every $m$-quasi-Einstein manifold satisfies
\begin{eqnarray}
\label{r}
\frac{1}{2}\Delta R-\frac{m+2}{2m}\langle\nabla f,\nabla R\rangle&=&-\frac{m-1}{m}\Big|Ric-\frac{R}{n}g\Big|^{2}\nonumber\\&&-\frac{n+m-1}{mn}(R-n\lambda)\Big(R-\frac{n(n-1)}{n+m-1}\lambda\Big).
\end{eqnarray} In particular, in the steady case we have $Ric=0,$ provided that $R=0$. Whence relation (\ref{eq6}) remounts to  Calabi  \cite{calabi} and Yau \cite{yau}.

Proceeding, it is well-known that such a metric is real analytic (cf. Proposition 2.8 in \cite{petersen}), hence the zeroes of the scalar curvature $R$ are isolated. Therefore, if $R\ge0$, but $R\neq0,$ we choose $p\in M^n$ such that $R(p)>0$ and a ball $B_{p}(r_{0})$ with radius $r_{0}>0$ such that $$\int_{B_{p}(r_{0})}uRd\sigma=mC_{0}$$ is a positive constant. Then we use the trace of (\ref{eq5}) to conclude that for all $r\geq r_{0}$

\begin{equation*}
mC_{0}=\int_{B_{p}(r_{0})}uRd\sigma\leq\int_{B_{p}(r)}uRd\sigma= m\int_{B_{p}(r)}\Delta ud\sigma.
\end{equation*} Next we invoke Stokes formula and (\ref{eqvol4}) to deduce
\begin{eqnarray*}
mC_{0}&\leq& m\int_{\partial B_{p}(r)}\frac{\partial u}{\partial\eta}ds \leq m\int_{\partial B_{p}(r)}|\nabla u| ds\nonumber\\
&\leq& m\sqrt{\frac{\mu}{m-1}}\cdot Area(\partial B_{p}(r)).
\end{eqnarray*} This implies that for $r\geq r_{0}$ we have
\begin{equation}
\label{k1}
Area(\partial B_{p}(r))\geq c > 0
\end{equation} for an uniform constant $c$. 

Finally, on integrating (\ref{k1}) from $r_{0}$ to $r$ we arrive at
$$Vol(B_{p}(r))\geq c(r-r_{0})\geq c_{0}\cdot r$$
for all $r\geq 2r_{0}.$ This finishes the proof of the theorem.
\end{proof}

\subsection{Proof of Theorem \ref{thmE}}

\begin{proof}
Firstly we notice that using the hypotheses on $f$ and $\mu$ in (\ref{eq8}) we deduce
\begin{equation}
\label{k2}
|\nabla u|^{2}\leq\frac{-\lambda e^{2k/m}}{m-1}.
\end{equation} Next, taking into account that $R\geq\lambda n$ and $u>0$ we infer $$\int_{B_{p}(r)} u(R-\lambda n)d\sigma \geq 0$$ for each $r>0.$ Moreover, if for all $r>0$ we have $$\int_{B_{p}(r)} u(R-\lambda n)d\sigma=0,$$ then $R=\lambda n$ on $M^n$. Hence from (\ref{r}) we deduce that $M^n$ is Einstein, but this gives a contradiction. In fact, according to \cite{Case} does not exist Einstein manifold with non-trivial expanding quasi-Einstein structure and potential function bounded from below.

From now on the proof look likes that one of the previous theorem. In particular, there is $p\in M^n$ such that $R(p)>\lambda n$. Since $u$ is positive there exists $r_{0}>0$ such that $$\int_{B_{p}(r_{0})}u(R-\lambda n)d\sigma=mC_{0}$$ is a positive constant.  Moreover, from the analyticity of $R$ it follows that for all $r\geq r_{0}$

\begin{eqnarray}
\label{poi}
mC_{0}&\leq&\int_{B_{p}(r)}u(R-\lambda n)d\sigma= m\int_{B{p}(r)}\Delta ud\sigma\nonumber\\
&=& m\int_{\partial B_{p}(r)}\frac{\partial u}{\partial\eta}ds \leq m\int_{\partial B_{p}(r)}|\nabla u|ds \nonumber\\
&\leq& m\sqrt{\frac{-\lambda e^{2k/m}}{m-1}}\cdot Area(\partial B_{p}(r)),
\end{eqnarray} where we have used Stokes formula and (\ref{k2}). Therefore, (\ref{poi}) allows us to deduce that for $r\geq r_{0}$
\begin{equation}
\label{k3}
Area(\partial B_{p}(r))\geq c > 0
\end{equation} for an uniform constant $c.$ 

In order to conclude it suffices to integrate (\ref{k3}) from $r_{0}$ to  $r$ to arrive at
$$Vol(B_{p}(r))\geq c(r-r_{0})\geq c_{0}\cdot r$$
for all $r\geq 2r_{0},$ which gives the requested result.
\end{proof}

\subsection{Proof of Theorem \ref{thmE-C}}
\begin{proof}
First, since $\mu=0$ we use (\ref{2eq}) to infer
\begin{equation}
\label{k4}
\Delta e^{-f}=(-\Delta f + |\nabla f|^{2})e^{-f}=-\lambda m e^{-f}.
\end{equation}
Now, upon integrating of (\ref{k4}) over $B_{p}(r)$ we deduce

\begin{eqnarray}
\label{k5}
-\lambda m\int_{B_{p}(r)}e^ {-f}d\sigma&=&\int_{B_{p}(r)}\Delta e^{-f}d\sigma\nonumber\\&=&\int_{\partial B_{p}(r)}\frac{\partial}{\partial \eta}(e^{-f})ds.
\end{eqnarray}
On the other hand, from \cite{wang} we have $|\frac{\partial f}{\partial \eta}|\leq |\nabla f|\leq\sqrt{-\lambda m}.$ This enables us to use (\ref{k5}) to arrive at
\begin{equation}
\label{eq999}
-\lambda m\int_{B_{p}(r)}e^ {-f}d\sigma\leq \sqrt{-\lambda m}\int_{\partial B_{p}(r)}e^{-f}ds.
\end{equation} Next, denoting
$$\xi(r):=Vol_{f}(B_{p}(r))=\int_{B_{p}(r)}e^{-f}d\sigma,$$
we can use (\ref{eq999}) to get
$$\xi'(r)\geq\sqrt{-\lambda m}\xi(r).$$
Whence, on integrating this inequality from $1$ to $r$ we conclude that $$\xi(r)=\int_{B_{p}(r)}e^{-f}d\sigma\geq ce^{\sqrt{-\lambda m}r}$$ for any $c>0.$ From here we conclude the proof of the theorem.
\end{proof}

\subsection{Proof of Theorem \ref{thmPW}}
\begin{proof}
We start invoking Theorem \ref{1teo}  to deduce that  $(M^{n},g,\nabla f,\lambda)$ is an expanding $m$-quasi-Einstein manifold with potential function $f=-m\ln u$ and $\mu<0.$

Next, we combine $f$-volume estimates obtained by Qian \cite{qian} and Theorem 9 in \cite{pigola} to conclude that the weak Maximum Principle at infinity is valid for the $f$-Laplacian on $(M^{n},g,\nabla f,\lambda).$ 
We also highlight that the potential function $f$ satisfies $|\nabla f|^{2}\leq -\lambda m$ (cf. \cite{wang}). From this setting, we apply the weak Maximum Principle at infinity for $|\nabla f|^{2}$ to conclude that there exists a sequence $\{p_{k}\}\subset M^{n},$ such that
\begin{equation*}
|\nabla f|^{2}(p_{k})\geq\overline{ |\nabla f|}^{2}-\frac{1}{k} \quad\quad and \quad\quad \Delta_{f}|\nabla f|^{2}(p_{k})\leq \frac{1}{k},
\end{equation*} where $\overline{|\nabla f|}^{2}=\sup_{M}|\nabla f|^{2}.$

We now recall the weighted Bochner formula:
\begin{equation}
\label{090}
\frac{1}{2}\Delta_{f}|\nabla f|^{2}=|Hess f|^{2}+Ric(\nabla f,\nabla f)+Hess f(\nabla f,\nabla f)+\langle\nabla f,\nabla\Delta_{f}f\rangle.
\end{equation} Therefore, (\ref{eqqem}) and (\ref{2eq}) substituted in (\ref{090}) gives
\begin{eqnarray}
\label{901}
\frac{1}{2}\Delta_{f}|\nabla f|^{2}&\geq &\lambda |\nabla f|^{2}+\frac{1}{m}|\nabla f|^{4}-2\mu e^{2f/m}|\nabla f|^{2}\nonumber\\
&\geq &\big(\lambda-2\mu u^{-2}+\frac{1}{m}|\nabla f|^{2}\big)|\nabla f|^{2}.
\end{eqnarray} Now, since $u\leq\sqrt{\frac{2\mu}{\lambda}}$ we immediately deduce
\begin{equation}
\frac{1}{2}\Delta_{f}|\nabla f|^{2}\geq\frac{1}{m}|\nabla f|^{4}.
\end{equation} From this, over $\{p_{k}\}$ we get
$$\frac{1}{2k}\geq \frac{1}{m}\Big(\overline{|\nabla f|}^{2}-\frac{1}{k}\Big)^{2}.$$
So, when $k$ goes to infinity we conclude that $\overline{|\nabla f|}^{2}=0$ and this forces $f$ to be constant which concludes the proof of the theorem.
\end{proof}

\begin{acknowledgement}
The E. Ribeiro Jr  would like to thank the Department of Mathematics - Lehigh University, where part of this work was carried out. He is grateful to Huai-Dong Cao for the warm hospitality and his constant encouragement. Moreover, he wish thank Lin Feng Wang for making available his article. Finally, the authors want to thank the referee for his careful reading and helpful suggestions. 
\end{acknowledgement}

\end{document}